\documentclass[11pt]{article}

\usepackage{graphicx}
\usepackage{color,verbatim}
\newtheorem{thm}{Theorem}[section]

\newtheorem{pro}[thm]{Proposition}
\newtheorem{lem}[thm]{Lemma}
\newtheorem{cor}[thm]{Corollary}
\newtheorem{alg}[thm]{Algorithm}
\newtheorem{ass}[thm]{Assumption}
\newtheorem{defi}[thm]{Definition}
\newtheorem{example}[thm]{Example}
\newcommand{\sect}[1]{
        \par
        \stepcounter{section}
        \settowidth{\hangindent}{\large\bf\thesection.~}
        \hangafter=1
        \bigskip\bigskip\noindent
        {\large\bf\hbox{\thesection.~}#1}\par
        \nopagebreak
        \medskip
        \renewcommand{\theequation}{\thesection.\arabic{equation}}
        \setcounter{equation}{0}
        \setcounter{subsection}{0}
}
\renewcommand{\subsection}[1]{
        \stepcounter{subsection}
        \noindent
        {\bf\hbox{\thesubsection.~}#1}
        \nobreak
}
\renewcommand{\subsubsection}[1]{
        \stepcounter{subsubsection}
        \noindent
        {\bf\hbox{\thesubsubsection.~}#1}
        \nobreak
}

\newcommand{\nn}{\nonumber}
\newcommand{\be}{\begin{equation}}
\newcommand{\ee}{\end{equation}}
\newcommand{\ba}{\begin{array}}
\newcommand{\ea}{\end{array}}
\newcommand{\bea}{\begin{eqnarray}}
\newcommand{\eea}{\end{eqnarray}}
\newcommand{\bal}{\begin{alg}}
\newcommand{\eal}{\end{alg}}
\newcommand{\ble}{\begin{lem}}
\newcommand{\ele}{\end{lem}}
\newcommand{\bco}{\begin{cor}}
\newcommand{\eco}{\end{cor}}
\newcommand{\bde}{\begin{defi}}
\newcommand{\ede}{\end{defi}}
\newcommand{\bth}{\begin{thm}}
\newcommand{\eth}{\end{thm}}
\newcommand{\bpr}{\begin{pro}}
\newcommand{\epr}{\end{pro}}
\newcommand{\bas}{\begin{ass}}
\newcommand{\eas}{\end{ass}}
\newcommand{\bex}{\begin{example}}
\newcommand{\eex}{\end{example}}
\newcommand{\reff}[1]{(\ref{#1})}
\newcommand{\refa}[1]{Assumption\ \ref{#1}}
\newcommand{\refl}[1]{Lemma\ \ref{#1}}
\newcommand{\reft}[1]{Theorem\ \ref{#1}}

\newcommand{\refc}[1]{Corollary\ \ref{#1}}

\newcommand{\refp}[1]{Proposition\ \ref{#1}}
\newcommand{\refal}[1]{Algorithm\ \ref{#1}}

\newtheorem{prc}[thm]{Procedure}
\newcommand{\bprc}{\begin{prc}}
\newcommand{\eprc}{\end{prc}}
\def\dd{&\!\!\!\!}
\def\eop{\hfill\vbox{\hrule height0.6pt\hbox{\vrule height1.3ex
width0.6pt\hskip1.2ex\vrule width0.6pt}\hrule height0.6pt}}
\def\prf{\noindent {\sl Proof.} \rm}
\def\alglist{
\begin{list}{Step 1}
{\setlength{\leftmargin}{0.5 in}\setlength{\labelwidth}{0.7 in}}
}
\def\eli{\end{list}}
\def\bdes{\begin{description}}
\def\edes{\end{description}}

\def\na{\nabla}

\def\la{\lambda}

\def\hf{\frac{1}{2}}

\def\alp{\alpha}

\def\st{\hbox{s.t.}}
\def\diag{\hbox{diag}\,}

\begin{document}
\pagenumbering{arabic}
\begin{titlepage}\setcounter{page}{0}

\title{A primal-dual interior-point relaxation method \\
for nonlinear programs
}
\author{
Xin-Wei Liu\thanks{Institute of Applied Mathematics, Hebei University of Technology, Tianjin 300401, China. E-mail:
mathlxw@hebut.edu.cn. The research is supported by the Chinese NSF grants (nos. 11671116 and 11271107), the Major Research Plan of the National Natural Science Foundation of China (no. 91630202), and a grant (no. A2015202365) from
Hebei Natural Science Foundation of China.} \
and Yu-Hong Dai\thanks{Academy of Mathematics and Systems Science, Chinese Academy of Sciences, Beijing 100190, China.
 This author is supported by the Chinese NSF grants (nos. 11631013, 11331012 and 81173633) and the National Key Basic Research Program of China (no.
2015CB856000).}
}
\date{ }
\maketitle

\noindent\underline{\hspace*{6.3in}}
\par

\vskip 10 true pt \noindent{\small{\bf Abstract.}
We prove that the classic logarithmic barrier problem is equivalent to a particular logarithmic barrier positive relaxation problem with barrier and scaling parameters.
Based on the equivalence, a line-search primal-dual interior-point relaxation method for nonlinear programs is presented.
Our method does not require any primal or dual iterates to be interior-points, which is prominently different from
the existing interior-point methods in the literature.
A new logarithmic barrier penalty function dependent on both primal and dual variables is used to prompt the global convergence of the method, where the penalty parameter is updated adaptively.
Without assuming any regularity condition, it is proved that our method will terminate at an approximate KKT point of
the original problem provided the barrier parameter tends zero. Otherwise, either an approximate infeasible stationary point or an approximate singular stationary point of the original problem will be found. Some preliminary numerical results are reported, including the results for a well-posed problem for which many line-search interior-point methods were demonstrated not to be globally convergent,
a feasible problem for which the LICQ and the MFCQ fail to hold at the solution and an infeasible problem, and for some standard test problems of the CUTE collection.
These results show that our algorithm is not only efficient for well-posed feasible problems, but also is applicable for some ill-posed feasible problems and some even infeasible problems.

\noindent{\bf Key words:} nonlinear programming, constrained optimization, interior-point method, logarithmic barrier problem,
global convergence

\noindent{\bf AMS subject classifications.} 90C26, 90C30, 90C51

\noindent\underline{\hspace*{6.3in}}

\vfil\eject
}
\end{titlepage}

\sect{Introduction}

We consider the general nonlinear programs with equality and inequality constraints
\bea \hbox{minimize}\quad(\min)\dd\dd f(x) \label{probo}\\
\hbox{subject to}\quad(\st)\dd\dd h_i(x)=0,\ i=1,\ldots,m_e, \label{probic}\\
                           \dd\dd c_j(x)\le 0,\ j=1,\ldots,m,\label{probec} \eea where $m_e$ and $m$ are integer numbers,
$x\in\Re^n$, $f$, $h_i (i=1,\ldots,m_e)$, and $c_j (j=1,\ldots,m)$ are twice continuously differentiable real-valued functions defined on $\Re^{n}$.
Interior-point methods have been among those most efficient methods for nonlinear programs (for example, see \cite{ByrHrN99,curtis12,GoOrTo03,VanSha99,WacBie06}). Given parameter $\mu>0$,
interior-point methods for nonlinear program \reff{probo}--\reff{probec} usually
need to approximately solve the following logarithmic barrier problem
\bea \min\dd\dd f(x)-\mu\sum_{j=1}^{m}\ln t_j \label{bo1}\\
\st\dd\dd h_i(x)=0,\ i=1,\ldots,m_e, \label{bic1}\\
\dd\dd c_j(x)+t_j= 0,\ j=1,\ldots,m, \label{bec1}
\eea
where $t_j (j=1,\ldots,m)$ are slack variables for inequality constraints.
Some other interior-point methods approximately solve the KKT system of the above logarithmic barrier problem as follows,
\bea
\dd\dd\na f(x)+\sum_{i=1}^{m_e}\lambda_i\na h_i(x)+\sum_{j=1}^{m}s_j\na c_j(x)= 0, \label{bkkt11}\\
\dd\dd t_j>0,\ s_j>0,\ t_js_j-\mu= 0,\ j=1,\ldots,m, \label{bkkt12}\\
\dd\dd h_i(x)= 0,\ i=1,\ldots,m_e;\quad c_j(x)+t_j=0,\ j=1,\ldots,m, \label{bkkt13}\eea
where $\la=(\la_i, i=1,\ldots,m_e)\in\Re^{m_e}$, $s=(s_j, j=1,\ldots,m)\in\Re^{m}$ are respectively the Lagrangian multiplier vector associated with constraints in \reff{bic1} and \reff{bec1}, $\mu>0$ is the barrier parameter which can be any number of a decreasing sequence with the limit $0$.

During the iterative processes for the logarithmic barrier problem \reff{bo1}--\reff{bec1} and the system \reff{bkkt11}--\reff{bkkt13}, all iterates must be interior points. That is, there should always have $t=(t_j, j=1,\ldots,m)>0$ and $s>0$. This guarantee for iterates to be interior-points may result in the truncation of the step and can impact on global performance of many existing line-search interior-point methods. An analytical counterexample presented in W\"achter and Biegler \cite{WacBie00} demonstrated that interior-point methods using linearized constraints and interior-point guarantee may fail in converging to any feasible point of a well-posed problem. As $\mu>0$ is small enough, an approximate KKT point of problem \reff{bo1}--\reff{bec1} is usually thought to be an approximate KKT point of the original problem.

There are already many interior-point methods in the literature which do not suffer from the failure in \cite{WacBie00}. These methods either use trust region techniques for new iterates (such as \cite{ByrHrN99}) or change the system \reff{bkkt11}--\reff{bkkt13} and/or its linearized system (e.g. \cite{curtis12,ForGil98,GoOrTo03,LiuSun01,LiuYua07,UlbUlV04,WacBie06}). Very recently, with help of the augmented Lagrangian function, Dai, Liu and Sun \cite{DLS17} presented a new primal-dual interior-point method for nonlinear programs. The method can produce interior-point iterates without truncation of the step, which is entirely different from the current primal-dual interior-point methods in the literature. It is proved that the method converges to a KKT point of the original problem as the barrier parameter tends to zero.
Otherwise, the penalty parameter tends to zero, and the method converges to either an infeasible stationary point or a singular stationary point of
the original problem. In particular, the method has the capability to rapidly detect the infeasibility of the problem. It is a very useful property in practice.

In this paper, we firstly prove that the logarithmic barrier problem \reff{bo1}--\reff{bec1} is equivalent to a particular logarithmic barrier positive relaxation problem with barrier and scaling parameters (see problem \reff{probo1}--\reff{probic2} in the next section). That is, we can derive a KKT point of problem \reff{bo1}--\reff{bec1} by solving its equivalent problem \reff{probo1}--\reff{probic2}.
Based on the equivalence, a primal-dual interior-point relaxation method for nonlinear programs \reff{probo}--\reff{probec} is then presented.
Our method does not require any primal or dual iterates to be interior-points, which is prominently different from the existing interior-point methods in the literature. This characteristic of our method makes us avert the truncation of any step for guarantee of interior-point iterates.

Our method has some similarity to that in \cite{DLS17}. However, since our method is proposed without using the augmented Lagrangian function, it has more flexibility to the selection of penalty parameter. Moveover, our method in this paper solves general optimization problems with inequality and equality constraints. It is known that the Lagrangian multipliers are dependent on the scaling of constraints. Thus, a scaling parameter is incorporated into the method. The incorporation of this parameter makes our method be robust for some degenerate and difficult problems.
Without assuming any regularity condition, it is proved that our method will terminate at an approximate KKT point of
the original problem provided the barrier parameter tends to zero. Otherwise, either an approximate infeasible stationary point or an approximate singular stationary point of the original problem will be found. Some preliminary numerical results are reported, including the results for  a well-posed problem for which many line-search interior-point methods were demonstrated not to be globally convergent, a feasible problem for which the LICQ and the MFCQ fail to hold at the solution and an infeasible problem, and for some standard test problems of the CUTE collection. These results show that our algorithm is not only efficient for well-posed feasible problems, but also is applicable for some ill-posed feasible problems and some even infeasible problems.

This article is organized as follows. In section 2, we describe a particular logarithmic barrier positive relaxation problem, and prove its equivalence to the logarithmic barrier problem \reff{bo1}--\reff{bec1}. Our primal-dual interior-point relaxation method for nonlinear program \reff{probo}--\reff{probec} is presented in section 3. In section 4, we show the global convergence results of our method. Some preliminary numerical results are reported in section 5. We conclude the paper in section 6.

Throughout the article, we use standard notations from the literature. A letter with
subscript $k$ ($l$) is related to the $k$th ($l$th)
iteration, the subscript $j$ ($i$) indicates the $j$th ($i$th) component of a vector, and the subscript $kj$ is the $j$th
component of a vector at the $k$th iteration. All vectors are column vectors, and $z=(x,u)$ means $z=[x^T,\hspace{2pt}u^T]^T$. The expression
$\theta_k={O}(t_k)$ means that there exists a constant $M$
independent of $k$ such that $|\theta_k|\le M|t_k|$ for all $k$ large enough, and
$\theta_k={o}(t_k)$ indicates that $|\theta_k|\le\epsilon_k|t_k|$ for all $k$ large enough with $\lim_{k\to
0}\epsilon_k=0$. If it is not specified, $I$ is an identity matrix whose order may be showed in the subscript or be clear in the context, $\|\cdot\|$ is the Euclidean norm. Some unspecified notations may be
identified from the context.

\sect{A logarithmic barrier positive relaxation problem}

Suppose that $\mu>0$ and $\tau>0$ are fixed constants. Let us consider problem
\reff{bo1}--\reff{bec1} and its KKT system
\reff{bkkt11}--\reff{bkkt13}. For $x\in\Re^n$, $t\in\Re^{m}$, and $s\in\Re^{m}$, define $z\in\Re^{m}$ and $y\in\Re^{m}$ by components \bea
z_j=(\sqrt{(\tau s_j-t_j)^2+4\tau\mu}-(\tau s_j-t_j))/2, \quad
y_j=(\sqrt{(\tau s_j-t_j)^2+4\tau\mu}+(\tau s_j-t_j))/2, \label{zydf}\eea
where $j=1,\ldots,m$. That is, both
$z:\Re^{2m}\to\Re^m$ and $y:\Re^{2m}\to\Re^m$ are functions on
$(t,s)$ and depend on the parameters $\mu$ and $\tau$. Obviously, $z_jy_j=\tau\mu$ for
$j=1,\ldots,m$. Moreover, we have following simple but important results on
$z$ and $y$. \ble\label{lemzp} For given $\mu>0$ and $\tau>0$, $z_j$ and $y_j$
are defined by \reff{zydf}. Then \bea t_j>0,\ s_j>0,\ t_js_j=\mu
\quad{\rm if\ and\ only\ if}\quad z_j-t_j=0,\ y_j-\tau s_j=0. \eea\ele
\prf Due to \reff{zydf}, $z_j>0$ and $y_j>0$ for any $t\in\Re^m$ and $s\in\Re^m$. If $z_j-t_j=0$, then $t_j>0$ and
$$\sqrt{(\tau s_j-t_j)^2+4\tau\mu}=t_j+\tau s_j.$$
Thus, $t_js_j=\mu$, which implies $s_j>0$. Similarly, one can prove that the results hold provided $y_j-\tau s_j=0$.

Conversely, suppose $t_j>0$ and $\ s_j>0$. If either $z_j-t_j\ne 0$ or $y_j-\tau s_j\ne 0$, then $t_js_j\ne\mu$, which shows that
$t_js_j=\mu$ implies $z_j-t_j=0$ and $y_j-\tau s_j=0$.
\eop

\ble\label{yp} For $j=1,\ldots,m$, let $z_j$ and $y_j$ be defined by \reff{zydf}. \\
(1) $z_j$ and $y_j$ are differentiable on $(t,s)$, and \bea
&&\na_tz_j=\frac{z_j}{z_j+y_j}e_j, \quad \na_ty_j=-\frac{y_j}{z_j+y_j}e_j, \label{20140327a}\\
&&\na_sz_j=-\tau\frac{z_j}{z_j+y_j}e_j, \quad
\na_sy_j=\tau\frac{y_j}{z_j+y_j}e_j, \label{20140327b} \eea
where $e_j\in\Re^m$ is the $j$-th coordinate vector.\\
(2) $z_j$ is a monotonically increasing function on $t_j$, and is a monotonically decreasing function on $s_j$.\\
(3) $y_j$ is a monotonically decreasing function on $t_j$, and is a monotonically increasing function on $s_j$.\\
(4) Both $z_j$ and $y_j$ are decreasing as $\mu>0$ becomes
smaller. \ele\prf
(1) Due to $z_jy_j=\tau\mu$, \bea
y_j\na_tz_j+z_j\na_ty_j=0. \label{lem1a1}\eea
Note that $z_j-y_j=t_j-\tau s_j$. Thus, $\na_t z_j-\na_t y_j=e_j$. Using \reff{lem1a1}, one has \reff{20140327a} immediately.

We can prove \reff{20140327b} similarly.

(2) By (1), $\frac{\partial z_j}{\partial t_j}>0$ and $\frac{\partial z_j}{\partial s_j}<0$.  The results can be derived straightforward.

(3) The results follows immediately since $\frac{\partial y_j}{\partial t_j}<0$ and $\frac{\partial y_j}{\partial s_j}>0$.

(4) Due to $z_jy_j=\tau\mu$ and $z_j-y_j=t_j-\tau s_j$, one has
\bea y_j\frac{dz_j}{d\mu}+z_j\frac{dy_j}{d\mu}=\tau,\ {\rm and}\  \frac{dz_j}{d\mu}-\frac{dy_j}{d\mu}=0. \nn\eea
Thus, \bea \frac{dz_j}{d\mu}=\frac{dy_j}{d\mu}=\frac{\tau}{z_j+y_j}>0, \label{lem2a1}\eea
which shows that both $z_j$ and $y_j$ are monotonically increasing on $\mu$.
\eop

Throughout this article, we denote $z:=z(t,s;\mu,\tau)$ and $y:=y(t,s;\mu,\tau)$ which are defined by
\reff{zydf} being functions on $(t,s)$. The next theorem provides a firm foundation for the development of our method.

\bth\label{mrs} Let $((x^*,t^*),(\la^*,s^*))$ be the KKT pair of problem
\reff{bo1}--\reff{bec1} and $(x^*,t^*,\la^*,s^*)$ satisfy the system
\reff{bkkt11}--\reff{bkkt13}, where $\la^*\in\Re^{m_e}$ and $s^*\in\Re^m$ are respectively the
Lagrangian multipliers of constraints \reff{bic1} and \reff{bec1}. Then
$((x^*,t^*,s^*),(\la^*,s^*,s^*))$ are the KKT pair of the following problem \bea
\min_{x,t,s}\dd\dd f(x)-\mu\sum_{j=1}^m \ln z_j \label{probo1}\\
\st\dd\dd h_i(x)=0,\ i=1,\ldots,m_e, \label{probic1}\\
   \dd\dd c_j(x)+t_j=0,\ j=1,\ldots,m, \label{probec1}\\
   \dd\dd z_j-t_j=0,\ j=1,\ldots,m.\label{probic2}\eea
That is, $\la^*\in\Re^{m_e}$, $s^*\in\Re^m$ and $s^*\in\Re^m$ are the corresponding Lagrangian multipliers of
constraints \reff{probic1}, \reff{probec1} and \reff{probic2}, respectively.

Conversely, if $((x^*,t^*,s^*),(\la^*,\beta^*,\nu^*))$ are the KKT pair of the
problem \reff{probo1}--\reff{probic2}, where $\la^*\in\Re^{m_e}$, $\beta^*\in\Re^m$ and
$\nu^*\in\Re^m$ are the associated Lagrangian multipliers of constraints
\reff{probic1}, \reff{probec1} and \reff{probic2}, then $\beta^*=\nu^*=y^*(t^*,s^*;\mu,\tau)/\tau=s^*$ and $(x^*,t^*,\la^*,s^*)$
satisfies the system \reff{bkkt11}--\reff{bkkt13}. Thus, $((x^*,t^*),(\la^*,s^*))$ is the KKT
pair of problem \reff{bo1}--\reff{bec1}. \eth\prf Due to
\refl{lemzp}, the KKT conditions of problem
\reff{bo1}--\reff{bec1} can be written as follows: \bea
\dd\dd\na f(x^*)+\sum_{i=1}^{m_e}\lambda_i^*\na h_i(x^*)+\sum_{j=1}^ms_j^*\na c_j(x^*)= 0, \label{bkkt31}\\
\dd\dd h_i(x^*)=0,\ i=1,\ldots,m_e;\quad c_j(x^*)+t_j^*=0,\ j=1,\ldots,m,\label{bkkt33}\\
\dd\dd z_j^*-t_j^*=0,\ j=1,\ldots,m,\label{bkkt34} \eea where $z^*=z(t^*,s^*;\mu,\tau)$ and $y^*=y(t^*,s^*;\mu,\tau)$.

Using \refl{yp}, we can derive the following KKT conditions of problem \reff{probo1}--\reff{probic2}: \bea
\dd\dd\na f(x^*)+\sum_{i=1}^{m_e}\lambda_i^*\na h_i(x^*)+\sum_{j=1}^m\beta_j^*\na c_j(x^*)=0, \label{bkkt41}\\
\dd\dd \beta^*-\sum_{j=1}^m\frac{\mu+y_j^*\nu_j^*}{z_j^*+y_j^*}e_j= 0, \label{bkkt41a}\\
\dd\dd \sum_{j=1}^m\tau\frac{\mu-z_j^*\nu_j^*}{z_j^*+y_j^*}e_j=0, \label{bkkt42}\\
\dd\dd h_i(x^*)=0,\ i=1,\ldots,m_e,\label{bkkt43}\\
\dd\dd c_j(x^*)+t_j^*=0,\ j=1,\ldots,m, \label{bkkt43a}\\
\dd\dd z_j^*-t_j^*=0,\ j=1,\ldots,m,\label{bkkt44} \eea
where $\la_i^* (i=1,\ldots,m_e)$, $\beta_j^* (j=1,\ldots,m)$ and $\nu_j^* (j=1,\ldots,m)$ are the Lagrangian multipliers of constraints \reff{probic1}, \reff{probec1} and \reff{probic2}, respectively.

For $j=1,\ldots,m$, $z_j^*\nu_j^*=\mu$ and $t_j^*s_j^*=\mu$ due to \reff{bkkt42} and \reff{bkkt44}. Thus, $\nu_j^*=s_j^*$ for all $j$ and $\beta^*=\nu^*$ by \reff{bkkt41a}.  Since $z_j^*y_j^*=\tau\mu$, $\nu_j^*=y_j^*/\tau$ for $j=1,\ldots,m$. Thus, any solution of the system \reff{bkkt41}--\reff{bkkt44} satisfies the system \reff{bkkt31}--\reff{bkkt34}.
Similarly, we can prove that any solution of the system \reff{bkkt31}--\reff{bkkt34} also solves the system \reff{bkkt41}--\reff{bkkt44}.
Hence, the result follows immediately from the equivalence between the systems \reff{bkkt31}--\reff{bkkt34} and \reff{bkkt41}--\reff{bkkt44}. \eop

Since $z$ is a function on $(t,s)$, problem \reff{probo1}--\reff{probic2} is a nonlinear programming problem on $(x,t,s)$ with equality constraints. Although the logarithmic barrier positive relaxation problem \reff{probo1}--\reff{probic2} has a similar form to the classic logarithmic barrier problem \reff{bo1}--\reff{bec1}, it is essentially distinguished from problem \reff{bo1}--\reff{bec1} in that the relaxation problem does not require $t$ and $s$ to be positive. This important characteristic makes us free to truncate the step for guarantee of interior-point iterates, a technique generally used by the existing interior-point methods in the literature.

\sect{A primal-dual interior-point relaxation algorithm}

Our algorithm consists of the inner algorithm and the outer algorithm. In the inner algorithm, we attempt to find an approximate KKT point of the logarithmic barrier positive relaxation problem \reff{probo1}--\reff{probic2} for any given $\mu>0$ and $\tau>0$. In the outer algorithm, we update parameters $\mu$ and $\tau$ and select the initial value of penalty parameter $\rho$ of the merit function according to the information of the solution derived from the inner algorithm.

Throughout the article, we denote $v:=(x,t,s)\in\Re^{n+2m}$ and functions $F:\Re^{n+2m}\to\Re$, $C:\Re^{n+2m}\to\Re^{m_e+2m}$,
$$F(v):=f(x)-\mu\sum_{j=1}^m \ln z_j,\ C(v):=(h(x),c(x)+t,z-t),$$
where we ignore the parameters $\mu$ and $\tau$ for simplicity of statement.
Then, due to \refl{yp},
\bea \dd\dd\na F(v)=\left(\na f(x),-\sum_{j=1}^m\frac{\mu}{z_j+y_j}e_j, \sum_{j=1}^m\frac{\tau\mu}{z_j+y_j}e_j\right), \nn\\
\dd\dd\na C(v)=\left(\ba{ccc}
\na h(x) & \na c(x) & 0 \\
0        & I_m      & -(Z+Y)^{-1}Y \\
0        & 0        & -\tau(Z+Y)^{-1}Z\ea\right), \nn
\eea
where $e_j\in\Re^m$ is the $j$-th coordinate vector in $\Re^m$, $I_m$ is the order-$m$ identity matrix, $Z=\diag(z)$ and $Y=\diag(y)$.

Let
$$L(v,w)=f(x)-\mu\sum_{j=1}^m \ln z_j+\sum_{i=1}^{m_e}\la_i h_i(x)+\sum_{j=1}^m\beta_j (c_j(x)+t_j)+\sum_{j=1}^m\nu_j(z_j-t_j)$$
be the Lagrangian function of problem \reff{probo1}--\reff{probic2}, where $w=(\la,\beta,\nu)\in\Re^{m_e+2m}$.
Note $z_jy_j=\tau\mu,\ j=1,\ldots,m$. By \refl{yp}, its Hessian has the form
\bea
\na_{vv}^2L(v,w)=\left[\ba{ccc}
H(x,\la,\beta) & 0 & 0 \\
0 & \mu\sum_{j=1}^m\frac{(\tau\nu_j+z_j)+(\tau\nu_j-y_j)}{(z_j+y_j)^3}e_je_j^T & -\tau\mu\sum_{j=1}^m\frac{(\tau\nu_j+z_j)+(\tau\nu_j-y_j)}{(z_j+y_j)^3}e_je_j^T \\[8pt]
0 &-\tau\mu\sum_{j=1}^m\frac{(\tau\nu_j+z_j)+(\tau\nu_j-y_j)}{(z_j+y_j)^3}e_je_j^T & \tau^2\mu\sum_{j=1}^m\frac{(\tau\nu_j+z_j)+(\tau\nu_j-y_j)}{(z_j+y_j)^3}e_je_j^T\ea\right], \nn\eea
where $H(x,\la,\beta)=\na^2 f(x)+\sum_{i=1}^{m_e}\la_i\na^2 h_i(x)+\sum_{j=1}^m\beta_j\na^2c_j(x)$. If we take $\nu=y/\tau$ in $w$, then
\bea \na_{vv}^2L(v,w)=\left[\ba{ccc}
H(x,\la,\beta) & 0 & 0 \\
0 & \sum_{j=1}^m\frac{\mu}{(z_j+y_j)^2}e_je_j^T & -\sum_{j=1}^m\frac{\tau\mu}{(z_j+y_j)^2}e_je_j^T \\[8pt]
0 & -\sum_{j=1}^m\frac{\tau\mu}{(z_j+y_j)^2}e_je_j^T & \sum_{j=1}^m\frac{\tau^2\mu}{(z_j+y_j)^2}e_je_j^T\ea\right]. \label{170819a}\eea

\subsection{The subproblems for search direction.}
Suppose $v_k:=(x_k,t_k,s_k)$ be the current iterate and $w_k:=(\la_k,\beta_k,\nu_k)$ be the corresponding estimate of the Lagrangian multiplier. The classic SQP approach for problem \reff{probo1}--\reff{probic2} solves the quadratic programming (QP) subproblem
\bea \min\dd\dd\na F(v_k)^Td+\hf d^TQ(v_k,w_k)d \label{qpf}\\
\st\dd\dd C(v_k)+\na C(v_k)^Td=0, \label{qpc} \eea
where $Q(v_k,w_k)$ is some positive definite approximation to the Hessian $\na_{vv}^2L(v_k,w_k)$.
Note that, even though $\na c(x_k)$ is of full column rank, the constraint system \reff{qpc} may still have a unbounded solution
since some diagonal element in $(Z_k+Y_k)^{-1}Z_k$ may be close to zero. In order to keep the solution $d_k$ bounded in our algorithm, we introduce the well-behaved null-space technology to the constraint system (for example, see \cite{byrd,ByrGiN00} for trust-region methods and \cite{LiuSun01,LiuYua00} for line-search methods). Moreover, $Q(v_k,w_k)$ is selected to have the same form as $\na_{vv}^2L(v_k,w_k)$ in \reff{170819a}, where $H(x_k,\la_k,\beta_k)$ is replaced by a positive definite approximation $B_k\in\Re^{n\times n}$ to $H(x_k,\la_k,\beta_k)$. Thus, $Q(v_k,w_k)$ is an approximation to $\na_{vv}^2L_(v_k,w_k)$ and
$$Q(v_k,w_k)=\left[\ba{ccc}
B_k & 0 & 0 \\
0 & \sum_{j=1}^m\frac{\mu}{(z_{kj}+y_{kj})^2}e_je_j^T & -\sum_{j=1}^m\frac{\tau\mu}{(z_{kj}+y_{kj})^2}e_je_j^T \\[8pt]
0 & -\sum_{j=1}^m\frac{\tau\mu}{(z_{kj}+y_{kj})^2}e_je_j^T & \sum_{j=1}^m\frac{\tau^2\mu}{(z_{kj}+y_{kj})^2}e_je_j^T\ea\right].$$

Motivated by above arguments, we firstly approximately solve the subproblem
\bea \min\dd\dd q_k^N(d;\rho):=\hf\rho d^TQ(v_k,w_k)d+\|C(v_k)+\na C(v_k)^Td\| \label{qsubo1}\\[8pt]
\st\dd\dd\|Rd\|\le\xi\|R^{-1}\na C(v_k)C(v_k)\|, \label{qsubc1}\eea
where $\xi>1$ is a constant, $\rho>0$ is the current value of the penalty parameter used in the merit function (see the next subsection), $R=\diag(1,\ldots,1,\tau,\ldots,\tau)\in\Re^{(n+2m)\times(n+2m)}$ with $1$'s of number $(n+m)$ and $\tau$'s of number $m$, respectively.
Let $p_k\in\Re^{n+2m}$ be the solution. Then our search direction $d_k$ is generated by the null-space QP subproblem
\bea \min\dd\dd q_k(d):=\na F(v_k)^Td+\hf d^TQ(v_k,w_k)d  \label{mqpf}\\
\st\dd\dd \na C(v_k)^T(d-p_k)=0. \label{mqpc}\eea
Apparently, $q_k(d_k)\le q_k(p_k)$. In particular, if $C(v_k)=0$, $p_k=0$ and $q_k(d_k)\le 0$.
\ble\label{lem0622} Assume $\na C(v_k)C(v_k)\ne 0$. Let $p_k$ be a solution of subproblem \reff{qsubo1}--\reff{qsubc1}. If $Q(v_k,w_k)$ is positive semi-definite, \bea \label{eq:2.17}
\dd\dd \|C(v_k)\|-q_k^N(p_k;\rho) \nn\\
\dd\dd\ge\frac{1}{2}\min\{1,\eta_k\}\frac{\|R^{-1}\na C(v_k)C(v_k)\|^2}{\|C(v_k)\|^2}\{\|C(v_k)\| \label{060501}\\
\dd\dd\hspace{10pt}-\rho\|C(v_k)\|^2\frac{C(v_k)^T\na C(v_k)^TR^{-2}Q(v_k,w_k)R^{-2}\na C(v_k)C(v_k)}{\|R^{-1}\na C(v_k)C(v_k)\|^2}\}, \nn
\eea where $\eta_k=\|R^{-1}\na C(v_k)C(v_k)\|^2/\|\na C(v_k)^TR^{-2}\na C(v_k)C(v_k)\|^2$. Thus,
\bea \|C(v_k)\|-q_k^N(p_k;\rho)\ge\frac{1}{4}\min\{1,\eta_k\}\frac{\|R^{-1}\na C(v_k)C_{\mu}(v_k)\|^2}{\|C(v_k)\|}\label{20170622b}\eea provided $\rho\|C(v_k)\|\le1/(2\lambda_{\max}(R^{-1}Q(v_k,w_k)R^{-1}))$, where $\lambda_{\max}(R^{-1}Q(v_k,w_k)R^{-1})$ is the maximal eigenvalue of $R^{-1}Q(v_k,w_k)R^{-1}$.
\ele\prf Since $p^C_k=-\min(1,\eta_k)R^{-2}\na C(v_k)C(v_k)$ is feasible to the subproblem \reff{qsubo1}--\reff{qsubc1},
\bea q_k^N(p_k;\rho)\le q_k^N(p_k^C;\rho). \nn\eea
Thus, $\|C(v_k)\|-q_k^N(p_k;\rho)\ge \|C(v_k)\|-q_k^N(p_k^C;\rho)$. Due to the positive semi-definiteness of $Q(v_k,w_k)$, \bea
\dd\dd\|C(v_k)\|-q_k^N(p_k;\rho) \nn\\
\dd\dd\ge\|C(v_k)\|-\|C(v_k)+\na C(v_k)^Tp_k^C\|\nn\\
\dd\dd\hspace{10pt}-\hf\rho\min\{1,\eta_k\}C(v_k)^T\na C(v_k)^TR^{-2}Q(v_k,w_k)R^{-2}\na C(v_k)C(v_k). \nn
\eea
By Lemma 2.1 of Liu and Yuan \cite{LiuYua07}, \bea \|C(v_k)\|-\|C(v_k)+\na C(v_k)^Tp_k^C\|\ge \frac{1}{2}\min\{1,\eta_k\}\|R^{-1}\na C(v_k)C(v_k)\|^2/\|C(v_k)\|. \nn\label{20170622a}\eea
The result \reff{eq:2.17} follows immediately.

If $2\rho\lambda_{\max}(R^{-1}Q(v_k,w_k)R^{-1})\|C(v_k)\|\le 1$, then
$$\rho\|C(v_k)\|^2\frac{C(v_k)^T\na C(v_k)^TR^{-2}Q(v_k,w_k)R^{-2}\na C(v_k)C(v_k)}{\|R^{-1}\na C(v_k)C(v_k)\|^2}\le\hf\|C(v_k)\|.$$
Thus, \reff{20170622b} is derived from \reff{eq:2.17}. \eop

\subsection{The merit function.} The merit function plays an important role in prompting the global convergence of the algorithm. For the logarithmic barrier positive relaxation problem \reff{probo1}--\reff{probic2}, we introduce the merit function
\bea \phi(v;\rho)=\rho (f(x)-\mu\sum_{j=1}^m\ln{z_j})+\|(h(x),c(x)+t,z-t)\|, \nn\eea
where $\rho>0$ is a penalty parameter, parameters $\mu$ and $\tau$ are ignored for simplicity of statement. It is a function on $(x,t,s)$ (where $x\in\Re^n$ and $t\in\Re^m$ are primal variables and $s\in\Re^m$ is a dual vector), thus it is essentially different from that used in some existing primal-dual interior-point methods such as \cite{CheGol06,LiuSun01,LiuYua07}. It is noted that some methods based on augmented Lagrangian have introduced merit functions with primal and dual variables, for example, see \cite{DLS17,ForGil98,GilRob12,wright95}. However, they generated their search directions by different systems and subproblems and shared different motivations with our algorithm. In view of previous notations in this section, it can also be written as
$$\phi(v;\rho)=\rho F(v)+\|C(v)\|.$$

\ble\label{mefunc} Given parameters $\mu>0$ and $\tau>0$. For any $\rho\ge 0$, $v\in\Re^{n+2m}$, $d\in\Re^{n+2m}$ and $d\ne 0$, the directional derivative $\phi^{'}(v,d;\rho)$ of function $\phi(v;\rho)$ at $v$ along $d$ exists, and \bea \phi^{'}(v,d;\rho)\le\pi(v,d;\rho), \label{ddphi}\eea
where $\pi(v,d;\rho)=\rho\na F(v)^Td+\chi(v,d)$ and $\chi(v,d)=\|C(v)+\na C(v)^Td\|-\|C(v)\|$.
\ele\prf
If $C(v)\ne 0$, $\phi(v;\rho)$ is differentiable at $v$. Thus, $\phi^{'}(v,d;\rho)$ exists for any $d\ne 0$.
For $v$ such that $C(v)=0$, using the definition of the directional derivative (for example, see (A.14) of \cite{NocWri99}), one has $\phi^{'}(v,d;\rho)=\rho\na F(v)^Td+\|\na C(v)^Td\|$. Therefore, $\phi^{'}(v,d;\rho)$ exists for every $v\in\Re^{n+2m}$, $d\in\Re^{n+2m}$ and $d\ne 0$.

Since \bea
\dd\dd\phi^{'}(v,d;\rho)=\lim_{\alp\downarrow 0}{\phi(v+\alp d;\rho)-\phi(v;\rho)\over\alp} \nn\\
\dd\dd=\rho\na F(v)^Td+\lim_{\alp\downarrow 0}{\|C(v)+\alp\na C(v)^Td+{o}(\alp)\|-\|C(v)\|\over\alp} \nn\\
\dd\dd\le\rho\na F(v)^Td+\|C(v)+\na C(v)^Td\|-\|C(v)\|+\lim_{\alp\downarrow 0}{\|{o}(\alp)\|\over \alp} \nn\\
\dd\dd=\rho\na F(v)^Td+\|C(v)+\na C(v)^Td\|-\|C(v)\|, \nn
\eea
\reff{ddphi} follows immediately. \eop

The following result shows that $d_k$ generated by subproblem \reff{mqpf}--\reff{mqpc} can be a descent direction of the merit function $\phi(v;\rho)$, provided $\rho$ is suitably selected.
\ble For given $\mu>0$ and $\tau>0$, suppose that $d_k$ is a solution of the QP subproblem \reff{mqpf}--\reff{mqpc} at $v_k$. Let $\phi^{'}_{k}(d_{k};\rho)=\phi^{'}(v_k,d_k;\rho)$, $\delta\in (0,1)$ is a constant. If $\rho$ is small enough, then $$\phi^{'}_{k}(d_{k};\rho)\le(1-\delta)(q_k^N(p_k;\rho)-\|C(v_k)\|)-\hf\rho d_k^TQ(v_k,w_k)d_k.$$ In particular,
$\phi^{'}_{k}(d_{k};\rho)\le -\hf\rho d_k^TQ(v_k,w_k)d_k$ for any $\rho>0$ provided $\|C(v_k)\|=0$.
\ele\prf
Due to \refl{mefunc}, $\phi^{'}_{k}(d_{k};\rho)\le\rho\na F(v_k)^Td_k+\chi(v_k,d_k)$. Thus, \bea
\phi^{'}_{k}(d_{k};\rho)\le\rho\na F(v_k)^Tp_k-\hf\rho d_k^TQ(v_k,w_k)d_k+q_k^N(p_k;\rho)-\|C(v_k)\| \label{20170622c} \eea
since $\na F(v_k)^Td_k+\hf d_k^TQ(v_k,w_k)d_k\le\na F(v_k)^Tp_k+\hf p_k^TQ(v_k,w_k)p_k$.
If $\|C(v_k)\|=0$, $p_k=0$ and $q_k^N(p_k;\rho)=0$. Thus, $\phi^{'}_{k}(d_{k};\rho)\le -\hf\rho d_k^TQ(v_k,w_k)d_k$.
Otherwise, by \refl{lem0622}, one can take $\rho$ small enough such that $q_k^N(p_k;\rho)-\|C(v_k)\|<0$ and $$\rho\na F(v_k)^Tp_k+\delta(q_k^N(p_k;\rho)-\|C(v_k)\|)\le 0.$$ Then the result follows immediately from \reff{20170622c}.
\eop

\ble\label{cypp} Given $\mu>0$ and $\tau>0$. For any $j=1,\ldots,m$, $z_{k+1,j}\ge t_{k+1,j}$
if either $t_{k+1,j}\le 0$ or $t_{k+1,j}>0$ but $s_{k+1,j}\le\frac{\mu}{t_{k+1,j}}$, where $z_{k+1,j}=z_j(t_{k+1},s_{k+1};\mu,\tau)$ is given by \reff{zydf}.
\ele\prf
If $t_{k+1,j}\le 0$, then $z_{k+1,j}-t_{k+1,j}>0$ since $z_{k+1,j}>0$. Now we consider the case $t_{k+1,j}>0$.
If $\tau s_{k+1,j}+t_{k+1,j}\le 0$, one has $s_{k+1,j}<0<\frac{\mu}{t_{k+1,j}}$; else $\tau s_{k+1,j}+t_{k+1,j}>0$, $z_{k+1,j}\ge t_{k+1,j}$ if and only if \bea
\sqrt{(\tau s_{k+1,j}-t_{k+1,j})^2+4\tau\mu}\ge\tau s_{k+1,j}+t_{k+1,j}, \nn\eea
which is equivalent to $t_{k+1,j}s_{k+1,j}\le\mu.$ The result follows immediately since $t_{k+1,j}>0$.
\eop

\subsection{Our algorithm.}
Similar to the existing primal-dual interior-point methods, our algorithm consists of the inner algorithm and the outer algorithm, where the inner algorithm tries to find an approximate KKT point of the logarithmic barrier problem, while the outer algorithm updates the parameters by the information derived from the inner algorithm.

For convenience of statement, we denote
\bea \dd\dd r(v_k,\la_k;\mu,\tau)=\left[\ba{c}
\na f(x_k)+\na h(x_k)\la_k+\na c(x_k)s_k \\
C(v_k)\ea\right],\nn\\[8pt]
\dd\dd g(v_k;\mu,\tau)=\frac{1}{\|C(v_k)\|}\left[\ba{c}
\na h(x_k)h(x_k)+\na c(x_k)(z_k-t_k) \\
c(x_k)+t_k-(z_k-t_k) \\
Z_k(z_k-t_k)\ea\right],\nn
\eea
where $v_k=(x_k,t_k,s_k)$, $z_k=z(t_k,s_k;\mu,\tau)$, $C(v_k)=(h(x_k),c(x_k)+t_k,z_k-t_k)$, $Z_k=\diag(z_k)$.
The following proposition provides the terminating conditions for our algorithm.
\bpr\label{prop1} For any given $\mu>0$ and $\tau>0$, let $\{v_k\}$ be any sequence convergent to the point $v^*$, where $v_k=(x_k,t_k,s_k)\in\Re^n\times\Re^m\times\Re^m$, $v^*=(x^*,t^*,s^*)$. Let $\epsilon>0$ be any given small scalar. \\
(1) If $(v^*,w^*)$ is a KKT pair of problem \reff{probo1}--\reff{probic2}, where $w^*=(\la^*,\beta^*,\nu^*)\in\Re^{m_e}\times\Re^m\times\Re^m$ is the corresponding Lagrangian multiplier, then there exists a $\la_k\in\Re^{m_e}$ such that $\|r(v_k,\la_k;\mu,\tau)\|_{\infty}\le\epsilon$ for every sufficiently large $k$. \\
(2) If $z_k\ge t_k$ for all $k\ge0$ and
\bea \lim_{k\to\infty}R^{-1}\na C(v_k)C(v_k)/\|C(v_k)\|=0, \label{171022a}\eea
then $\|g(v_k;\mu,\tau)\|_{\infty}\le\epsilon$ for every sufficiently large $k$. \\
\epr\prf
(1) If $(v^*,w^*)$ is a KKT pair of problem \reff{probo1}--\reff{probic2}, $(v^*,w^*)$ satisfies the KKT conditions \reff{bkkt41}--\reff{bkkt44}. Thus, $\beta^*=\nu^*=s^*$, and $t^*>0$, $t^*_js_j^*=\mu$ for $j=1,\ldots,m$. It follows from \reff{bkkt41} and \reff{bkkt43a}, \reff{bkkt44} that $x^*$, $t^*$, $\la^*$, and $s^*$ satisfy \reff{bkkt11}--\reff{bkkt13}. Hence, $\|r_1(v^*,\la^*;\mu,\tau)\|_{\infty}=0$. The result follows immediately from the continuity of $\na f(x)$, $\na h(x)$, $\na c(x)$, and $C(v;\mu,\tau)$.

(2) Due to \reff{171022a}, one has
\bea \dd\dd \lim_{k\to\infty}\na h(x_k)h(x_k)/\|C(v_k)\|+\na c(x_k)(c(x_k)+t_k)/\|C(v_k)\|=0, \label{171029a}\\
     \dd\dd \lim_{k\to\infty}(c(x_k)+t_k)/\|C(v_k)\|-(Z_k+Y_k)^{-1}Y_k(z_k-t_k)/\|C(v_k)\|=0, \label{171029b}\\
     \dd\dd \lim_{k\to\infty}-(Z_k+Y_k)^{-1}Z_k(z_k-t_k)/\|C(v_k)\|=0. \label{171029c}\eea
By \reff{171029b} and \reff{171029c}, $\lim_{k\to\infty}(c(x_k)+t_k)/\|C(v_k)\|-(z_k-t_k)/\|C(v_k)\|=0$. Thus, \reff{171029a} implies
$$\lim_{k\to\infty}\frac{1}{\|C(v_k)\|}[\na h(x_k)h(x_k)+\na c(x_k)(z_k-t_k)]=0.$$
Furthermore, if \reff{171029c} holds, $\lim_{k\to\infty}Z_k(z_k-t_k)/\|C(v_k)\|=0$. Hence, the result is proved.
\eop

Now we are ready to present our primal-dual interior-point relaxation algorithm for problem \reff{probo}--\reff{probec}.
\bal\label{alg1} (A primal-dual interior-point relaxation algorithm for problem \reff{probo}--\reff{probec}) \ {\small \alglist
\item[Step 1] Given $(x_0,t_0,s_0)\in\Re^{n}\times\Re^m\times\Re^m$, $B_0\in\Re^{n\times n}$, $\mu_0>0$, $\tau_0>0$, $\rho_0>0$, $\delta\in(0,1)$,
$\sigma\in(0,\frac{1}{2})$, $\epsilon>0$.
Set $l:=0$.

\item[Step 2] While $\mu_l>\epsilon$ and $\tau_l>\epsilon$, start the inner algorithm. Otherwise, stop the algorithm.

Step 2.0 Let $(x_0,t_0,s_0)=(x_l,t_l,s_l)$, $B_0=B_l$, $\rho_0=\rho_l$, $\mu=\mu_l$ and $\tau=\tau_l$. Evaluate $z_0$ and $y_0$ by

\hspace{1cm} \reff{zydf}. Set $k:=0$.

Step 2.1 Find first an approximate solution $p_k$ of subproblem \reff{qsubo1}--\reff{qsubc1} such that \reff{eq:2.17} holds,

\hspace{1cm} then solve the QP subproblem \reff{mqpf}--\reff{mqpc} to derive $(d_{xk},d_{tk},d_{sk})$.

Step 2.2 {Choose} $\rho_{k+1}$ with either $\rho_{k+1}=\rho_k$ or $\rho_{k+1}\le0.5\rho_k$ such that
\bea &2\rho_{k+1}\|C(v_k)\|\frac{C(v_k)^T\na C(v_k)^TR^{-2}Q(v_k,w_k)R^{-2}\na C(v_k)C(v_k)}{\|R^{-1}\na C(v_k)C(v_k)\|^2}\le 1,\ {\rm and}&\label{rhoup0}\\[6pt]
&\pi(v_{k},d_{k};\rho_{k+1})\le(1-\delta)(q_k^N(p_k;\rho_{k+1})-\|C(v_k)\|)-\hf\rho_{k+1}d_k^TQ(v_k,w_k)d_k,\quad&
\label{rhoup}\eea
\hspace{1cm}where $\pi(v_{k},d_{k};\rho_{k+1})$ is defined in \refl{mefunc}.

Step 2.3 {Choose} the step-size $\alp_k\in (0,1]$ to be the maximal in $\{1,\delta,\delta^2,\ldots\}$ such that
\bea
\phi(x_k+\alp_kd_{xk},t_k+\alp_kd_{tk},s_k+\alp_kd_{sk};\rho_{k+1})-\phi(x_k,t_k,s_k;\rho_{k+1})\le\sigma\alpha_k\pi(v_{k},d_{k};\rho_{k+1}).
\label{srule}\eea

Step 2.4 {Set}  $x_{k+1}=x_k+\alpha_kd_{xk}$, $t_{k+1}=t_k+\alpha_kd_{tk}$, and $\hat s_{k+1}=s_k+\alp_kd_{sk}$.

Step 2.5 For $j=1,\ldots,m$, set \bea s_{k+1,j}=\left\{\ba{ll}
\hat s_{k+1,j}, & {\rm if}\quad t_{k+1,j}\le 0; \\[5pt]
\min\{\hat s_{k+1,j},\frac{\mu}{t_{k+1,j}}\}, & {\rm otherwise}. \ea\right.\label{uupdate}\eea

Step 2.6 Compute an estimate $\la_{k+1}$ of Lagrangian multiplier vector corresponding to the equality

\hspace{1cm} constraints in \reff{probic}.

Step 2.7 If $\|r(v_{k+1},\la_{k+1};\mu_l,\tau_l)\|_{\infty}\le 10\mu_l$, set $\mu_{l+1}=\min\{0.5\mu_l,\|r_1(v_{k+1},\la_{k+1};\mu_l,\tau_l)\|_{\infty}^{1.8}\}$,

\hspace{1cm} $\tau_{l+1}=\tau_l$.

\hspace{1cm} Else if $\|g(v_{k+1};\mu_l,\tau_l)\|_{\infty}\le\tau_l$, set $\mu_{l+1}=\mu_l$, $\tau_{l+1}\le0.6\tau_l$.

\hspace{1cm} In above two cases, stop the inner algorithm, set $(x_{l+1},t_{l+1},s_{l+1})=(x_{k+1},t_{k+1},s_{k+1})$,

\hspace{1cm} $l:=l+1$.

\hspace{1cm} Otherwise, update $B_k$ to $B_{k+1}$, set $k:=k+1$ and go to Step 2.1.

End (while)

\eli} \eal

Due to \refl{cypp}, the update \reff{uupdate} guarantees $z_{k+1}-t_{k+1}\ge 0$ for all $k$. Moreover, it will further reduce the value of $\phi(x_k+\alp_kd_{xk},t_k+\alp_kd_{tk},s_k+\alp_kd_{sk};\rho_{k+1})$. Thus, one has
\bea \phi(x_{k+1},t_{k+1},s_{k+1};\rho_{k+1})-\phi(x_k,t_k,s_k;\rho_{k+1})\le\sigma\alpha_k\pi(v_{k},d_{k};\rho_{k+1})<0
\label{srulea}\eea
for all $k\ge0$.

\sect{Global convergence}

Without assuming any regularity of constraints such as feasibility or constraint qualification, the local solution of the problem \reff{probo}--\reff{probec} may be either a KKT point or a singular stationary point (i.e. a Fritz-John point). For an infeasible problem, we usually want to find an infeasible stationary point which is a stationary point for minimizing some measure of constraint violations (see \cite{BurCuW14,ByrCuN10,NocOzW12}).

In order to solve the original problem \reff{probo}--\reff{probec}, \refal{alg1} approximately solves a sequence of logarithmic barrier positive relaxation problem \reff{probo1}--\reff{probic2}. In this section, we first prove that, for any given $\mu>0$ and $\tau>0$, the inner algorithm of \refal{alg1} will terminate in a finite number of iterations. Thus, either $\mu_l\to 0$ or $\tau_l\to 0$ as $l\to\infty$. After that, we consider the global convergence of the whole algorithm. It is proved that our algorithm will terminate at an approximate KKT point of the original problem provided the scaling parameter $\tau_l$ is far from zero. Otherwise, either an approximate infeasible stationary point or an approximate singular stationary point of the original problem will be found.

\subsection{Global convergence of the inner algorithm.} We consider the global convergence of the inner algorithm. Suppose that, for some given $\mu>0$ and $\tau>0$, the inner algorithm of \refal{alg1} does not terminate in a finite number of iterations and $\{(x_k,t_k,s_k)\}$ is an infinite sequence generated by the algorithm. We need the following blanket assumptions for our global convergence analysis.
\bas\label{ass1}\ \\
(1) The functions $f$ and $c_i (i\in{\cal I})$ are twice continuously
differentiable on $\Re^n$; \\
(2) The iterative sequence $\{x_k\}$ is in an open bounded set;\\
(3) The sequence $\{B_k\}$ is bounded, and for all $k\ge 0$ and $d_x\in\Re^n$, $d_x^TB_kd_x\ge\gamma\|d_x\|^2,$ where $\gamma>0$ is a constant;\\
(4) For all $k\ge 0$, $p_k$ is an approximate solution of subproblem \reff{qsubo1}--\reff{qsubc1} satisfying \reff{eq:2.17}.
\eas

We have the following results which are proved similar to Lemma 5 of \cite{ByrGiN00} and Lemma 4.2 of \cite{LiuSun01}.
\ble\label{lemabc} Suppose that \refa{ass1} holds. Then $\{z_k\}$ and $\{t_k\}$ are bounded, $\{y_k\}$ is componentwise bounded away from zero and $\{s_k\}$ is lower bounded. Furthermore, if the penalty parameter $\rho_k$ is bounded away from zero as $k\to\infty$, then $\{z_k\}$ is componentwise bounded away from zero, $\{y_k\}$ and $\{s_k\}$ are bounded. \ele\prf
\refa{ass1} implies that there exists a scalar $\chi>0$ such that $\|f(x_k)\|\le\chi$ and $\|c(x_{k})\|\le\chi$ for all $k\ge 0$.
Due to \reff{srulea}, $\phi(v_{k+1};\rho_{k+1})\le\phi(v_k;\rho_{k+1})$ for all $k\ge0$. Thus, for every $k\ge0$,
\bea \phi(v_{k+1};\rho_{k+1})-\phi(v_{k};\rho_{k})\le(\rho_k-\rho_{k+1})(\chi+m\mu\ln\|z_k\|). \nn\eea
Therefore, \bea
\phi(v_{k+1};\rho_{k+1})\le\phi(v_{0};\rho_{0})+(\rho_0-\rho_{k+1})(\chi+m\mu\max_{0\le l\le{k+1}}\ln\|z_l\|). \label{171018a}\eea
Note that the inequalities
\bea\phi(v_{k+1};\rho_{k+1})\ge-\rho_{k+1}(\chi+m\mu\max_{0\le l\le{k+1}}\ln\|z_l\|)+\|t_{k+1}\|-\|c(x_{k+1})\| \label{171018c}\eea and
$\phi(v_{k+1};\rho_{k+1})\ge-\rho_{k+1}(\chi+m\mu\max_{0\le l\le{k+1}}\ln\|z_l\|)+\|z_{k+1}\|-\|t_{k+1}\|$,
one has \bea
\phi(v_{k+1};\rho_{k+1})\ge-\rho_{k+1}(\chi+m\mu\max_{0\le l\le{k+1}}\ln\|z_l\|)+\|z_{k+1}\|-\|c(x_{k+1})\|. \label{171018b}\eea
Hence, it follows from \reff{171018a}, \reff{171018c}, and \reff{171018b} that, for all $k\ge 0$, \bea
\phi(v_{0};\rho_{0})+(1+\rho_{0})\chi+\rho_0m\mu\max_{0\le l\le{k+1}}\ln\|z_l\|)\ge\max(\|z_{k+1}\|,\|t_{k+1}\|), \nn\eea
which implies that $\{z_k\}$ is bounded. Furthermore, $\{t_k\}$ is bounded since $\{z_k\}$ is bounded.
Due to $z_{kj}y_{kj}=\tau\mu$ for every $j=1,\ldots,m$, the results on $\{y_k\}$ follow immediately.

For given $\mu>0$ and $\tau>0$, if $\{z_k\}$ is bounded, by \reff{zydf}, $s_{kj}>-\infty$ for all $k\ge 0$ and $j=1,\ldots,m$. Otherwise, if $s_{kj}\to-\infty$ for some $j$, then $z_{kj}\to\infty$, which is a contradiction. If $\{z_k\}$ is componentwise bounded away from zero, again by \reff{zydf}, $s_{kj}<+\infty$ for all $k\ge 0$ and $j=1,\ldots,m$. Thus, the results on $\{s_k\}$ are proved.
\eop

\bco\label{co1} Suppose that \refa{ass1} holds. Then, for any given $\mu>0$ and $\tau>0$, all sequences $\{C(v_k)\}$, $\{\na F(v_k)\}$, $\{\na C(v_k)\}$ and $\{Q(v_k,w_k)\}$ are bounded. Thus, there is a constant $\chi_0>0$ such that $\|C(v_k)\|\le\chi_0$, $\|\na C(v_k)\|\le\chi_0$, and $\lambda_{\max}(Q(v_k,w_k))\le\chi_0$.
\eco\prf For $j=1,\ldots,m$, due to $z_{kj}y_{kj}=\tau\mu$, $$\frac{1}{z_{kj}+y_{kj}}=\frac{z_{kj}}{z_{kj}^2+z_{kj}y_{kj}}\le\frac{1}{\tau\mu}z_{kj}.$$
Thus, by \refl{lemabc}, ${1}/({z_{kj}+y_{kj}})$ is bounded. The boundednesses of $\{C(v_k)\}$, $\{\na F(v_k)\}$, $\{\na C(v_k)\}$ and $\{Q(v_k,w_k)\}$ follow immediately from their expressions in previous section. \eop

\ble\label{lemg2} Suppose that \refa{ass1} holds. If
\bea {\|R^{-1}\na C(v_k)C(v_k)\|}\ge\chi_1{\|C(v_k)\|} \label{consc}
\eea
for some constant $\chi_1>0$ and for all $k\ge 0$, then there is a constant $\hat\rho>0$ such that $\rho_{k+1}=\hat\rho$ for all sufficiently large $k$.
\ele\prf Due to \refc{co1}, $\rho_{k+1}\|C(v_k)\|\le1/(2\lambda_{\max}(R^{-1}Q(v_k,w_k))R^{-1})$ provided $\rho_{k+1}\le\tau^2/(2\chi_0^2)$. Thus, if $\hat\rho\le\tau^2/(2\chi_0^2)$, it follows from \reff{20170622b} and \reff{consc} that \bea \|C(v_k)\|-\|C(v_k)+\na C(v_k)^Tp_k\|
\ge\frac{1}{4}\min\{1,\frac{1}{\chi_0^2}\}\chi_1^2{\|C(v_k)\|}. \label{20170703a}\eea

We achieve the result by proving that \reff{rhoup} holds with $\rho_{k+1}\le\hat\rho$ for some scalar $\hat\rho>0$. Since $q_k(d_k)\le q_k(p_k)$, one has
\bea &&\pi(v_{k},d_{k};\rho_{k+1})-(1-\delta)(q_k^N(p_k;\rho_{k+1})-\|C(v_k)\|)+\hf\rho_{k+1}d_k^TQ(v_k,w_k)d_k \nn\\
&&\le\rho_{k+1}(\na F(v_k)^Tp_k+\hf p_k^TQ(v_k,w_k)p_k)+\delta(\|C(v_k)+\na C(v_k)^Tp_k\|-\|C(v_k)\|)\nn\\
&&\le(\rho_{k+1}\xi_1-\delta\xi_2)\|C(v_k)\|,\nn \eea
where $\xi_1$ and $\xi_2$ are positive constants. Hence, \reff{rhoup} holds with $\rho_{k+1}\le\hat\rho$ provided  $\hat\rho=\min\{\tau^2/(2\chi_0^2),\delta\xi_2/\xi_1\}$.
\eop

\ble\label{imp} Suppose that \refa{ass1} holds, $d_k=(d_{xk},d_{tk},d_{sk})\in\Re^{n+2m}$ is a solution of QP \reff{mqpf}--\reff{mqpc},
$u_k=(\la_k,\beta_k,\nu_k)$ is the associated Lagrangian multiplier vector. If \reff{consc} holds for all sufficiently large $k$, the sequence $\{\|d_{k}\|\}$ is bounded. \ele\prf
If \reff{consc} holds for all sufficiently large $k$, by \refl{lemg2}, $\rho_{k}$ is bounded away from zero. It follows from \refl{lemabc} that $z_k$ and $y_k$ are bounded. Thus, there exists a constant $\chi_2>0$ such that \bea
\frac{\mu}{(z_{kj}+y_{kj})^2}\ge\chi_2\quad {\rm for}\quad j=1,\ldots,m. \nn \eea
Due to \bea
d^TQ(v_k,w_k)d&=& d_x^TB_kd_x+\sum_{j=1}^n\frac{\mu}{(z_{kj}+y_{kj})^2}({d_{tj}-\tau d_{sj}})^2 \nn\\
 &\ge& \gamma\|d_x\|^2+\chi_2\|d_t-\tau d_s\|^2 \label{Qpd}\eea  for every $d=(d_x,d_t,d_s)$, $d_k^TQ(v_k,w_k)d_k\ge\xi_3\|(d_{xk},d_{tk}-\tau d_{sk})\|^2$ for some scalar $\xi_3>0$.

Since $q_k(d_{k})\le q_k(p_k)$, $q_k(p_k)\le\chi_3$ and $q_k(d_{k})\ge-\chi_3\|(d_{xk},d_{tk}-\tau d_{sk})\|+\xi_3\|(d_{xk},d_{tk}-\tau d_{sk})\|^2$ for some scalars $\chi_3>0$ and $\xi_3>0$, one can deduce that $\|(d_{xk},d_{tk}-\tau d_{sk})\|$ is bounded. Otherwise, if $\|(d_{xk},d_{tk}-\tau d_{sk})\|$ is unbounded, then $\xi_3\le 0$, which is a contradiction. Thus, $\|d_{tk}\|$ is bounded due to $d_{tk}=p_{tk}-\na c(x_k)^T(d_{xk}-p_{xk})$. It implies that $\|d_{sk}\|$ is bounded. \eop

\ble\label{lemg3} Suppose that \refa{ass1} holds, $\{\alp_k\}$ is the sequence of step-sizes derived from \reff{srule} of \refal{alg1}. If the inequality \reff{consc} holds for all sufficiently large $k$, then $\{\alpha_k\}$ is bounded away from zero.
\ele\prf Due to Lemmas \ref{lemabc} and \ref{imp}, for every $j=1,\ldots,m$, one has
\bea &-\ln z_{j}(v_k+\alp d_k;\mu,\tau)+\ln z_{kj}-\alp\frac{1}{z_{kj}+y_{kj}}e_j^T(d_{tk}-\tau d_{sk})=o(\alp),& \\
&\|C(v_k+\alp d_k)\|=\|C(v_k)+\alp\na C(v_k)^Td_k\|+o(\alp)&
\eea
for all $\alp>0$ sufficiently small. Hence,  \bea\phi(v_k+\alp
d_{k};\rho_{k+1})-\phi(v_k;\rho_{k+1})=\alp\pi(v_k,d_{k};\rho_{k+1})+o(\alp)
\label{070121a}\eea for all $\alp\in[0,\tilde\alp]$, where $\tilde\alp>0$ is a sufficiently small constant. Due to
\bea
(1-\sigma)\alp\pi(v_k,d_{k};\rho_{k+1})\le
\alp(1-\sigma)(1-\delta)(q_k^n(p_k;\rho_{k+1})-\|C(v_k)\|)\le -\alp\xi_4\|C(v_k)\| \label{070116a}
\eea
(where $\xi_4=\xi_2(1-\sigma)(1-\delta)$), it follows from \reff{070121a} and \reff{070116a} that there
exists a scalar $\hat\alp\in (0,\tilde\alp]$ such that \bea \phi(v_k+\alp
d_{k};\rho_{k+1})-\phi(v_k;\rho_{k+1})\le\sigma\alp\pi(v_k,d_{k};\rho_{k+1})
\nn\eea for all $\alp\in(0,\hat\alp]$ and all $k\ge0$. Thus, by Step 2.3 of \refal{alg1}, $\alp_k\ge\delta\hat\alp$ for all $k\ge0$.
\eop

\ble\label{lemg4} Suppose that \refa{ass1} holds. If the inequality \reff{consc} holds for all sufficiently large $k$, then
\bea \lim_{k\to\infty}\|C(v_k)\|=0\quad{\rm and}\quad \lim_{k\to\infty}\|d_{k}\|=0. \eea
\ele\prf
According to \refl{lemg2}, without loss of generality, we suppose that $\rho_{k}=\hat\rho$ for all $k\ge 0$. Then, by \reff{srulea}, $\{\phi(v_k;\rho_{k+1})\}$ is a monotonically non-increasing sequence. Note that it is also a bounded sequence. Thus,
\bea \lim_{k\to\infty}\pi(v_k,d_{k};\rho_{k+1})=0 \eea
since $\alp_k$ is bounded away from zero. Using the last inequality of \reff{070116a}, one has
$$\lim_{k\to\infty}\|C(v_k)\|=0,$$ which implies $\lim_{k\to\infty}\|p_{k}\|=0$. Thus, by \reff{rhoup} and \reff{Qpd}, $\lim_{k\to\infty}\|d_{k}\|=0.$
\eop

Now we are ready to present our global convergence results on the inner algorithm.
\bth\label{lemj1} Given $\mu>0$ and $\tau>0$. Suppose that \refa{ass1} holds. Let $\{v_{k}\}$ and $\{\la_k\}$ be two sequences generated by the inner algorithm of \refal{alg1}. Let $\epsilon>0$ be any given small scalar. Then there is an integer $k>0$ such that either $\|r(v_{k+1},\la_{k+1};\mu,\tau)\|_{\infty}\le\epsilon$ or $\|g(v_{k+1};\mu,\tau)\|_{\infty}\le\epsilon$.
\eth\prf
If \reff{consc} holds for all $k\ge 0$, by \refl{lemg2}, $\rho_k$ remains a positive constant after a finite number of iterations. Thus, $\{z_k\}$ and $\{y_k\}$ are bounded above and componentwise bounded away from zero, $\{s_k\}$ is bounded.

We prove the result by contradiction. If the result does not hold, then the inner algorithm will not terminate in a finite number of iterations. Thus, $\{v_k\}$ is an infinite sequence.
Let $v^*=(x^*,t^*,s^*)$ be any limit point of $\{v_{k}\}$. Without loss of generality, suppose that $\lim_{k\to\infty} v_k=v^*$.
Due to \refl{lemg4}, $\lim_{k\to\infty}\|C(v_k)\|=0$ and $\lim_{k\to\infty}\|d_{k}\|=0$. Thus,
$$\lim_{k\to\infty}z_k=t^*>0,\quad \lim_{k\to\infty}y_k=\tau s^*>0$$
since $z_k-t_k=y_k-\tau s_k$. Moreover, by taking the limit on $k\to\infty$ in both sides of the following KKT condition
$$\na F(v_k)+Q(v_k,w_k)d_k+\na C(v_k)u_k=0$$
of subproblem \reff{mqpf}--\reff{mqpc}, due to \refc{co1}, one has
\bea \na F(v_k)+\na C(v_k)u_k=0. \label{170707a}
\eea
Note that $\lim_{k\to\infty}\|C(v_k)\|=0$. Therefore, every limit point of $\{v_k\}$ is a KKT point of \reff{probo1}--\reff{probic2}. In view of \refp{prop1}, there exists a $\la_k\in\Re^{m_e}$ such that the inequality $\|r(v_{k+1},\la_{k+1};\mu,\tau)\|_{\infty}\le\epsilon$ holds for every sufficiently large $k$.

In the following, suppose that \reff{consc} does not hold for all sufficiently large $k$. Then there exists some infinite index subset ${\cal K}$ such that the condition \reff{consc} does not hold for all $k\in{\cal K}$. That is,
\bea \lim_{k\in{\cal K},k\to\infty}{\|R^{-1}\na C(v_k)C(v_k)\|}/{\|C(v_k)\|}=0. \eea
Due to \refp{prop1}, one has $\|g(v_{k+1};\mu,\tau)\|_{\infty}\le\epsilon$ for all sufficiently large $k\in{\cal K}$.

The above arguments show that the sequence $\{v_k\}$ cannot be an infinite sequence. This contradiction implies the result of the theorem.
\eop

\subsection{Convergence results of the whole algorithm.} Now we consider the global convergence of the whole algorithm. It is well known that, without assuming any constraint qualification, a local solution of general nonlinear program can be either a KKT point or a Fritz-John point of the problem. For those nonlinear programs arising from practical situation, whether they are feasible are not known before solving them. Thus, some robust methods for nonlinear programs not only focus on convergence to KKT points of problems under some assumptions on constraint regularities, but also concern about convergence to Fritz-John points and infeasible stationary points of problems without assuming any regularity on constraints (for example, see \cite{BurCuW14,byrd,ByrCuN10,ByrMaN01,CheGol06,DLS17,LiuSun01,LiuYua07,yuan95}).

\bde\label{def2}
$x^*\in\Re^n$ is called a Fritz-John point or a singular stationary point of problem \reff{probo}--\reff{probec} if there exist $\la^*\in\Re^{m_e}$ and $\beta^*\in\Re^m$ such that
\bea && \na h(x^*)\la^*+\na c_i(x^*)\beta^*=0, \nn\\
&& h_i(x^*)=0,\ i=1,\cdots,m_e, \nn\\
&& \beta_j^*\ge 0,\ c_j(x^*)\le 0,\ \beta_j^*c_j(x^*)=0,\ j=1,\ldots,m. \nn\eea
\ede
\bde\label{def1}
$x^*\in\Re^n$ is called an infeasible stationary point of problem \reff{probo}--\reff{probec} if $x^*$ is an infeasible point and
\bea
\na h(x^*)h(x^*)+\na c(x^*)\max(0,c(x^*)=0. \nn
\eea\ede

Now we are ready to present our convergence results on the whole algorithm.
\bth\label{lemj2} Suppose that \refa{ass1} holds for every given parameters $\mu_l>0$ and $\tau_l>0$, sequences $\{x_l\}$ and $\{B_l\}$ are bounded. Let $\epsilon>0$ be small enough. Then either $\mu_l\le\epsilon$ or $\tau_l\le\epsilon$ for some $l$, and one of the following statement is true.  \\
(1) The parameter $\tau_l>\epsilon$, the inner algorithm terminates at a point $v_{l+1}$ where the terminating condition $\|r(v_{l+1},\la_{l+1};\mu_l,\tau_l)\|_{\infty}\le10\epsilon$ holds. \refal{alg1} terminates at an approximate KKT point of the original problem \reff{probo}--\reff{probec}; \\
(2) The parameter $\tau_l\le\epsilon$, the inner algorithm terminates at a point $v_{l+1}$ at which the condition $\|g(v_{l+1};\mu_l,\tau_l)\|_{\infty}\le\epsilon$ is satisfied and $\|C(v_k)\|$ is small enough. \refal{alg1} terminates at a point which is an approximate singular stationary point of the problem \reff{probo}--\reff{probec}; \\
(3) The parameter $\tau_l\le\epsilon$, the inner algorithm terminates due to $\|g(v_{l+1};\mu,\tau)\|_{\infty}\le\epsilon$ and $\|C(v_k)\|$ is bounded away from zero. \refal{alg1} terminates at a point $v_{l+1}$ which is an approximate infeasible stationary point of the problem \reff{probo}--\reff{probec}.
\eth\prf For every $l\ge 0$, \reft{lemj1} shows that, in \refal{alg1}, either $\mu_l$ or $\tau_l$ will be reduced. Finally, there is either $\mu_l\le\epsilon$ or $\tau_l\le\epsilon$. Furthermore,
The argument of \refl{lemabc} shows that $\{z_l\}$ and $\{t_l\}$ are bounded.

The condition $\|r(v_{l+1},\la_{l+1};\mu_l,\tau_l)\|_{\infty}\le\epsilon$ implies that $(x_{l+1},t_{l+1})$ is an approximate KKT point of
the logarithmic barrier positive relaxation problem \reff{bo1}--\reff{bec1}. It is typically known that, as $\mu_l$ is small enough, $x_{l+1}$ is also an approximate KKT point of the original problem \reff{probo}--\reff{probec}.

If it is other than the above case, then one has $\tau_{l}\le\epsilon$ and $\|g(v_{l+1};\mu_l,\tau_l)\|_{\infty}\le\epsilon$. That is,
\bea &&\|\na h(x_{l+1})h(x_{l+1})+\na c(x_{l+1})(z_{l+1}-t_{l+1})\|_{\infty}/\|C(v_{l+1})\|\le\epsilon, \label{170713a}\\
     &&\|Z_{l+1}(z_{l+1}-t_{l+1})\|_{\infty}/\|C(v_{l+1})\|\le\epsilon. \label{170713b}\eea
In order to demonstrate that the point $x_{l+1}$ satisfying \reff{170713a}--\reff{170713b}
is an approximate singular stationary point of the problem \reff{probo}--\reff{probec} when $\|C(v_{l+1})\|$ is small enough, we take the limit $\epsilon\to 0$ on both sides of \reff{170713a}--\reff{170713b}. The limit $\epsilon\to 0$ implies $\tau_l\to 0$.
Without loss of generality, suppose $z_{l+1}\to z^*$ and $x_{l+1}\to x^*$, $t_{l+1}\to t^*$ as $\tau_l\to 0$, where $x^*$, $t^*$, $z^*$ satisfy $h(x^*)=0$, $c(x^*)+t^*=0$, $z^*-t^*=0$. Furthermore, $h(x_{l+1})/\|C(v_{l+1})\|\to\la^*$, $(z_{l+1}-t_{l+1})/\|C(v_{l+1})\|\to\beta^*$,
Thus, $c(x^*)\le 0$, $\beta^*\ge 0$ (since $z_{l+1}-x_{l+1}\ge 0$ for all $l$),
and due to \reff{170713a}, \reff{170713b}, and $c(x^*)=-t^*=-z^*$,
\bea &\na h(x^*)\la^*+\na c(x^*)\beta^*=0,& \nn\\
     &\beta^*_jc_j(x^*)=0.&\nn\eea
That is, $x^*$ is a Fritz-John point problem \reff{probo}--\reff{probec}. Therefore, $x_{l+1}$ is an approximate singular stationary point of the problem.

If $\|g(v_{l+1};\mu_l,\tau_l)\|_{\infty}\le\epsilon$ but $C(v_{l+1})$ is bounded away from zero, then
\bea &&\|\na h(x_{l+1})h(x_{l+1})+\na c(x_{l+1})(z_{l+1}-t_{l+1})\|_{\infty}\le\epsilon, \label{171029a5}\\
     &&\|c(x_{l+1})+t_{l+1}-(z_{l+1}-t_{l+1})\|_{\infty}\le\epsilon, \label{171029c5}\\
     &&\|Z_{l+1}(z_{l+1}-t_{l+1})\|_{\infty}\le\epsilon. \label{171029b5}\eea
Similar to the preceding arguments, suppose $z_{l+1}\to z^*$ and $x_{l+1}\to x^*$, $t_{l+1}\to t^*$ as $\tau_l\to 0$, then
$t^*=\hf (z^*-c(x^*))$ due to \reff{171029c5}. Thus, by \reff{171029b5}, $z^*_j(z_j^*+c_j(x^*))=0$ for $j=1,\ldots,m$. This fact implies
\bea c_j(x^*)+t_j^*=z_j^*-t_j^*=\max(0, c_j(x^*), \quad j=1,\ldots,m. \nn\eea
Therefore, it follows from \reff{171029a} that
$$\na h(x^*)h(x^*)+\na c(x^*)\max(0,c(x^*))=0. $$
That is, $x^*$ is a stationary point to minimize $\hf\|(h(x),\max\{0,c(x)\})\|^2$, thus is an infeasible stationary point of problem \reff{probo}--\reff{probec}. It shows that $x_{l+1}$ is an approximate infeasible stationary point of problem \reff{probo}--\reff{probec}.\eop

\sect{Numerical experiments}

The numerical experiments were conducted on a Lenovo laptop with the LINUX operating system (Fedora 11).
The algorithm were implemented in MATLAB (version R2008a). Two kinds of test problems originated from the literature were solved, including some simple but hard problems, which may be an infeasible problem, a problem feasible but LICQ and MFCQ failing
to hold at the solution, or is a well-posed one but some class of interior-point methods was proved not to be globally convergent, and some standard test problems of the CUTE collection \cite{BonCGT95}. These problems have been solved in \cite{DLS17}.

We use the standard initial point $x_0$ for all test problems, and set $t_0=-c(x_0)$, $s_{0j}=\min\{1,0.95\mu/\max(0,t_{0j})\}$. The initial parameters are selected as follows: $\mu_0=0.1$, $\tau_0=1$, $\delta=0.5$, $\sigma=0.0001$, and $\epsilon=10^{-8}$. The initial
penalty is selected to be dependent on the initial point $x_0$ as $\rho_0=\min\{100,\max(1,\|(\max(0,c(x_0)),h(x_0))\|/|f(x_0)|)\}$, $B_0$ is simply taken as the identity matrix and $B_k$ is updated by the well-known Powell's damped BFGS update formula (for example, see \cite{NocWri99}).
The subproblems are solved by similar techniques as those used in \cite{DLS17}, while subproblem \reff{qsubo1}--\reff{qsubc1} in this paper should be treated more carefully since $Q(v_k,w_k)$ is always positive semi-definite. The
whole algorithm is terminated as either $\mu_l\le\epsilon$ or $\tau_l\le\epsilon$, or the total number of iterations (that is, the
number of solving QP \reff{mqpf}--\reff{mqpc}) is larger than 1000.

\subsection{Numerical results on three simple but hard problems.} In this subsection, we report our numerical results on three simple but hard examples taken from the literature. In Tables 1, 2 and 3, the number in column $l$ means that the data in the row are taken from the corresponding outer iterate at which either $\mu_l$ or $\tau_l$ is reduced, $f_l=f(x_l)$, $v_l=\|(h(x_l),\max(0,c(x_l))\|$, $\|r_l\|_{\infty}=\|r(x_{l+1},s_{l+1},\la_{l+1};\mu_l,\tau_l)\|_{\infty}$, $\|g_l\|_{\infty}=\|g(x_{l+1},s_{l+1};\mu_l,\tau_l)\|_{\infty}$, $k$ is the number of inner iterations needed from $(\mu_{l-1},\tau_{l-1})$ to $(\mu_l,\tau_l)$.

The first example was presented by W\"{a}chter and Biegler \cite{WacBie00} and further discussed by Byrd, Marazzi and Nocedal \cite{ByrMaN01}:
\begin{eqnarray}
\hbox{min} && x\sb{1} \nonumber\\
({\rm TP1}) \quad\quad
\hbox{s.t.} && x\sb{1}\sp{2}-x\sb{2}-1=0, \nonumber\\
   && x\sb{1}-x\sb{3}-2=0, \nonumber\\
   && x\sb{2}\geq 0, \quad x\sb{3}\geq 0. \nonumber
\end{eqnarray}
The standard initial point is $x_0=(-4,1,1)$. This problem is a well-posed problem. It has a unique global minimizer $(2,3,0)$, at which gradients of the active constraints are linearly independent, and MFCQ holds. However, \cite{WacBie00} showed that many line search interior-point methods maight fail to find the solution.

Our algorithm terminates at the approximate solution $x^*=(2.0000,3.0000,0.0000)$ together with $t^*=(3.0000, 0.0000)$ and $s^*=(0.0000, 1.0000)$ in totally $19$ iterations. The numbers of function and gradient evaluations are $20$ and $20$, respectively. See Table \ref{tab1} for more details on iterations, from there one can observe the superlinear convergence of $v_l$ and $\|r_l\|_{\infty}$.
\begin{table}
{\small
\begin{center}
\caption{Output for test problem (TP1)}\label{tab1} \vskip 0.2cm
\begin{tabular}{|c|c|c|c|c|c|c|c|}
\hline
$l$ & $f_l$ & $v_l$ & $\|r_l\|_{\infty}$ & $\|g_l\|_{\infty}$ & $\mu_l$ & $\tau_l$ & $k$  \\
\hline
\hline
0 & -4 & 14 & 14 & 7.6026 & 0.1000 & 1 & - \\
1 & -1.5235 & 3.5562 & 3.5477 & 0.9974 & 0.1000 & 0.6000 &  2 \\
2 & -1.2344 & 2.2938 & 0.9737 & 0.5926 & 0.1000 & 0.3600 &  2 \\
3 & -1.0298 & 2.2408 & 0.4768 & 0.3133 & 0.1000 & 0.2160 &  1 \\
4 & -0.4257 & 1.9679 & 0.1957 & 0.6123 & 0.0500 & 0.2160 & 5 \\
5 & 2.0245 & 1.6851e-04 & 0.0049 & 4.4642 & 0.0011 & 0.2160 & 6 \\
6 & 2.0011 & 5.4789e-04 & 1.1834e-04 & 3.9978 & 1.3479e-06 & 0.2160 &  1 \\
7 & 2.0000 & 1.1711e-06 & 2.5296e-07 & 4.0000 & 1.0000e-09 & 0.2160 & 1 \\
8 & 2.0000 & 1.8137e-12 & 5.6019e-11 & 4.0000 & - & - & 1 \\
\hline
\end{tabular}
\end{center}}
\end{table}

The second example is a standard test problem taken from \cite[Problem 13]{HocSch81}:
\begin{eqnarray}
\hbox{min} && (x\sb{1}-2)\sp{2}+x\sb{2}\sp{2} \nonumber\\
({\rm TP2}) \quad\quad
\hbox{s.t.} && (1-x\sb{1})\sp{3}-x\sb{2}\geq 0, \nonumber\\
   && x\sb{1}\geq 0, \quad x\sb{2}\geq 0. \nonumber
\end{eqnarray}
The standard initial point $x_0=(-2,-2)$ for problem (TP2) is an infeasible point.
Note that its optimal solution $x^*=(1,0)$ is not a KKT point but a singular stationary point, at which LICQ and MFCQ fail to hold.

This problem has not been solved in \cite{ShaVan00,yamash98}, but has been solved in \cite{ByrHrN99,tseng99}.
\refal{alg1} terminated at an approximate point to the solution $x_l=(0.9905,-0.0000)$. The associated vectors are $t_l=(0.0000,0.9905,0.0000)$,
$s_l=1e3(7.3986,0.0000,7.3989)$. The numbers of iterations, function evaluations and gradient evaluations are $28$,
$76$ and $29$, respectively (see Table \ref{tab2}).

\begin{table}
{\small
\begin{center}
\caption{Output for test problem (TP2)}\label{tab2} \vskip 0.2cm
\begin{tabular}{|c|c|c|c|c|c|c|c|}
\hline
$l$ & $f_l$ & $v_l$ & $\|r_l\|_{\infty}$ & $\|g_l\|_{\infty}$ & $\mu_l$ & $\tau_l$ & $k$  \\
\hline
\hline
0 & 20 & 2 & 8.9116 & 0.7071 & 0.1000 & 1 & - \\
1 & 19.6742 & 2.7923 & 1.6715 & 0.7592 & 0.1000 & 1.5502e-05 &  1 \\
2 & 3.9492 & 0 & 6.2765e-06 & 1.0000 & 1.7127e-07 & 1.5502e-05 &  9 \\
3 & 1.0192 & 9.1401e-11 & 1.3450e-11 & 1.0000 & 1.0000e-09 & 1.5502e-05 &  17 \\
4 & 1.0192 & 9.1390e-11 & 1.3450e-11 & 0.0203 & - & - & 1 \\
\hline
\end{tabular}
\end{center}}
\end{table}

The third example is an infeasible problem named {\sl{isolated}} presented by Byrd, Curtis and Nocedal \cite{ByrCuN10}.
\bea
\min\dd\dd x\sb{1}+x_2 \nonumber\\
({\rm TP3}) \quad\quad
\st\dd\dd x_{1}^{2}-x\sb{2}+1\le 0, \nonumber\\
   \dd\dd x_{1}^{2}+x\sb{2}+1\le 0, \nonumber\\
   \dd\dd -x_{1}+x\sb{2}^2+1\le 0, \nonumber\\
   \dd\dd x_{1}+x\sb{2}^2+1\le 0. \nonumber
\eea
The standard initial point is $x_0=(3,2)$, its solution $x^*=(0,0)$ is a strict minimizer of the Euclidean norm of constraint infeasibility measures. The algorithm presented in \cite{ByrCuN10} found this point. Our algorithm terminates at an approximate point to it with $x_l=1.0e-04(-0.1547,-0.6259)$, $t_l=(-0.5000,-0.5000,-0.5000,-0.5000)$, $s_l=1.0e04(5.3752,-5.4001,3.6816,-3.6977)$ in $17$ iterations. Note that $t_l$ is not positive as it should always be in existing interior-point methods based on solving logarithmic barrier problem \reff{bo1}--\reff{bec1}.  The numbers of function and gradient evaluations are $20$ and $18$, respectively. Some more details please refer to Table \ref{tab3}.
\begin{table}
{\small
\begin{center}
\caption{Output for test problem (TP3)}\label{tab3} \vskip 0.2cm
\begin{tabular}{|c|c|c|c|c|c|c|c|}
\hline
$l$ & $f_l$ & $v_l$ & $\|r_l\|_{\infty}$ & $\|g_l\|_{\infty}$ & $\mu_l$ & $\tau_l$ & $k$  \\
\hline
\hline
0 & 5 & 12 & 43.1944 & 7.5805 & 0.1000 & 1 & - \\
1 & 2.5509 & 5.8735 & 4.3125 & 3.1566 & 0.1000 & 0.6000 & 1 \\
2 & 0.7102 & 2.3657 & 0.5663 & 0.5946 & 0.0500 & 0.6000 & 1 \\
3 & -0.1702 & 2.0217 & 0.4846 & 0.3160 & 0.0500 & 0.3600 & 1  \\
4 & 0.0606 & 2.0029 & 0.2179 & 0.1233 & 0.0500 & 0.2160 & 1  \\
5 & -0.0633 & 2.0042 & 0.1209 & 0.0903 & 0.0500 & 2.0359e-04 & 1  \\
6 & 1.0010e-04 & 2.0000 & 1.0183e-04 & 1.8054e-04 & 0.0500 & 1.0000e-09 & 11 \\
7 & -7.8060e-05 & 2.0000 & 5.0005e-10 & 2.3597e-07 & - & -  & 1 \\
\hline
\end{tabular}
\end{center}}
\end{table}

\subsection{Numerical results on test problems of the CUTE collection.}
A set of $59$ small- and medium-size test problems with general
inequality constraints from the CUTE collection \cite{BonCGT95}
were solved. These problems were selected since they had actual
numbers of problem variables and general inequality constraints
(i.e., not only bound constraints). Besides general inequality
constraints, some test problems may also have equality
constraints. In particular, these problems have been solved and
the numerical results have been reported in \cite{DLS17}.

We report our numerical results in Tables \ref{tab11} and
\ref{tab12}, where the columns ``$n$" and ``$m$" are the numbers
of variables and constraints of test problems, respectively. The columns of $f$ and $v$ show, respectively, the values of
objective functions and the infinite norms of constraint
violations at the terminating points, ``iter" represents the
total number of iterations needed for obtaining those values.
For comparison, we denote ``$\|r\|_{\infty}$" to be the same as ``$\|\phi\|_{\infty}$"
in \cite{DLS17}. The last two
columns ``$N_f$" and ``$N_g$" of Tables \ref{tab11} and
\ref{tab12} are respectively the numbers of evaluations of
functions and gradients needed by the algorithm.

Note that Tables \ref{tab11} and \ref{tab12} only list $48$ test problems for which the terminating conditions
of our algorithm have been satisfied before reaching the restriction of number of iterations.
The algorithm in \cite{DLS17} has terminated far from reaching
the restriction on the total number of iterations for $55$ problems.
The observed reason is that our algorithm in this paper seems to be more sensitive to the update of
scaling parameter corresponding to the penalty parameter in \cite{DLS17}, which possibly results in either more iterations
or smaller step-sizes. This may be reasonable when comparing with the interior-point method using an augmented Lagrangian function.

The preliminary results show that our algorithm can still be efficient
for most of test problems.
In particular, our algorithm need obviously fewer iterations, fewer evaluations of functions and gradients for some test problems
such as CHACONN2, HAIFAS, HS29, HS43, POLAK5, and so on.
However, our algorithm is still not available for using the exact Hessian,
thus it is not comparable with some well performed algorithms (for example, \cite{curtis12,WacBie06}) and software such as IPOPT.
Since our MATLAB implementation uses
MATLAB routines simply, it is believed that further improvements can
be achieved by using advanced techniques for, e.g., the computations
of subproblems.
\begin{table}
{\small
\begin{center}
\caption{Results for CUTE problems using approximate Hessians, part 1.}\label{tab11} \vskip 0.2cm
\begin{tabular}{|c|c|c|c|c|c|c|c|c|}
\hline
Problem & $n$ & $m$ & $f$ & $v$ & iter & $\|r\|_{\infty}$ & $N_f$ & $N_g$ \\
\hline
CB2 & 3 & 3 & 1.9522 & 9.4943e-08 & 8 & 8.2734e-06 & 9 & 9 \\
CB3 & 3 & 3 & 2.0000 & 0 & 9 & 4.3026e-09 & 10 & 10 \\
CHACONN1 & 3 & 3 & 1.9522 & 0 & 8 & 5.3363e-10 & 9 & 9 \\
CHACONN2 & 3 & 3 & 2.0000 & 0 & 9 & 2.6625e-11 & 10 & 10 \\
CONGIGMZ & 3 & 5 & 28.0000 & 2.5079e-09 & 36 & 3.4740e-07 & 69 & 37 \\
DIPIGRI & 7 & 4 & 680.6301 & 0 & 24 & 3.0065e-09 & 88 & 25 \\
EXPFITA & 5 & 22 & 0.0011 & 0 & 191 & 8.8791e-04 & 1478 & 192 \\
EXPFITC & 5 & 502 & 0.0233 & 0 & 707 & 0.0019 & 5813 & 708 \\
GIGOMEZ1 & 3 & 3 & -3.0000 & 0 & 16 & 6.5600e-11 & 43 & 17 \\
GIGOMEZ2 & 3 & 3 & 1.9522 &  9.7055e-10 & 8 & 9.6530e-08 & 10 & 9 \\
GIGOMEZ3 & 3 & 3 & 2.0000 & 0 & 13 & 7.5786e-08 & 28 & 14 \\
HAIFAS & 13 & 9 & -0.4500 & 0 & 12 & 1.7810e-09 & 21 & 13 \\
HS10 & 2 & 1 & -1.0000 & 3.6877e-09 & 11 & 8.0982e-09 & 12 & 12 \\
HS11 & 2 & 1 & -8.4985 & 0 & 7 & 3.8354e-09 & 8 & 8 \\
HS12 & 2 & 1 & -30.0000 & 0 & 12 & 8.0016e-09 & 36 & 13 \\
HS14 & 2 & 2 & 1.3935 & 0 & 6 & 4.5752e-13 & 7 & 7 \\
HS22 & 2 & 2 & 1.0000 & 0 & 8 & 6.4723e-11 & 10 & 9 \\
HS29 & 3 & 1 & -22.6274 & 0 & 35 & 2.0747e-14 & 142 & 36 \\
HS43 & 4 & 3 & -44.0000 & 0 & 14 & 1.1305e-09 & 18 & 15 \\
HS88 & 2 & 1 & 1.2980e-05 & 0.1332 & 6 & 2.0113e-10 & 8 & 7 \\
HS89 & 3 & 1 & 0.8163 & 0.0050 & 14 & 3.2978e-11 & 19 & 15 \\
HS90 & 4 & 1 & 7.8447e-06 & 0.1332 & 9 & 6.6620e-11 & 40 & 10 \\
HS91 & 5 & 1 & 8.2135e-05 & 0.1332 & 6 & 6.6644e-11 & 8 & 7 \\
HS92 & 6 & 1 & 2.6369e-06 & 0.1332 & 6 & 4.3302e-10 & 8 & 7 \\
HS100 & 7 & 4 & 680.6301 & 0 & 24 & 3.0065e-09 & 88 & 25 \\
HS100MOD & 7 & 4 & 678.6798 & 0 & 54 & 0.0036 & 263 & 55 \\
HS113 & 10 & 8 & 24.3062 & 0 & 139 & 2.3607e-10 & 1222 & 140 \\
KISSING & 127 & 903 & 32.7317 & 2.3054 & 85 & 1.4267e-08 & 1043 & 86 \\
KIWCRESC & 3 & 2 & 2.0000e-09 & 0 & 13 & 2.7207e-09 & 27 & 14 \\
MAKELA1 & 3 & 2 & -1.4142 & 0 & 10 & 1.4611e-09 & 22 & 11 \\
MAKELA2 & 3 & 3 & 7.2000 & 0 & 16 & 1.3945e-09 & 41 & 17 \\
MAKELA3 & 21 & 20 & 1.9895e-08 & 0 & 331 & 4.9710e-10 & 365 & 332 \\
MAKELA4 & 21 & 40 & 5.4155e-08 & 0 & 58 & 2.5122e-08 & 354 & 59 \\
MIFFLIN1 & 3 & 2 & -1.0000 & 0 & 7 & 4.0317e-07 & 8 & 8 \\
MINMAXRB & 3 & 4 & 2.1057e-06 & 6.0793e-06 & 508 & 0.0052 & 5492 & 509 \\
PENTAGON & 6 & 15 & 1.3685e-04 & 0 & 177 & 1.0524e-04 & 1066 & 178 \\
POLAK1 & 3 & 2 & 2.7183 & 0 & 165 & 1.6286e-09 & 1569 & 166 \\
POLAK3 & 12 & 10 & 5.9330 & 0 & 22 & 1.9982e-08 & 26 & 23 \\
POLAK5 & 3 & 2 & 50.0000 & 0 & 7 & 9.4837e-12 & 9 & 8 \\
POLAK6 & 5 & 4 & 15.6949 & 0.0013 & 52 & 6.6387e-10 & 553 & 53 \\
ROSENMMX & 5 & 4 & -44.0000 & 0 & 85 & 6.6579e-08 & 854 & 86 \\
S268 & 5 & 5 & 0.0289 & 0 & 21 & 3.2134e-05 & 106 & 22 \\
SPIRAL & 3 & 2 & -2.6781e-08 & 2.6781e-08 & 162 & 7.6149e-09 & 422 & 163 \\
\hline
\end{tabular}
\end{center}}
\end{table}
\begin{table}
{\small
\begin{center}
\caption{Results for CUTE problems using approximate Hessians, part 2.}\label{tab12} \vskip 0.2cm
\begin{tabular}{|c|c|c|c|c|c|c|c|c|}
\hline
Problem & $n$ & $m$ & $f$ & $v$ & iter & $\|r\|_{\infty}$ & $N_f$ & $N_g$ \\
\hline
TFI2 & 3 & 101 & 0.6468 & 0.0088 & 30 & 8.8354e-12 & 237 & 31 \\
VANDERM1 & 100 & 199 & 0 & 0.0168 & 6 & 1.0561e-10 & 8 & 7 \\
VANDERM2 & 100 & 199 & 0 & 0.0168 & 6 & 1.0561e-10 & 8 & 7 \\
VANDERM3 & 100 & 199 & 0 & 0.0203 & 7 & 1.6996e-10 & 8 & 8 \\
WOMFLET & 3 & 3 & 4.3957 & 0 & 23 & 1.3885e-07 & 116 & 24 \\
\hline
\end{tabular}
\end{center}}
\end{table}

\sect{Conclusion}

We present a primal-dual interior-point relaxation method for nonlinear programs in this article. The method is based on an equivalence between the classic logarithmic barrier problem and a particular logarithmic barrier positive relaxation problem. Remarkably different from the current primal-dual interior-point methods in the literature, our method does not require any primal or dual iterates to be interior-point points. Thus, the interior-point restriction in existing line-search methods can be removed. A new logarithmic barrier penalty function dependent on both primal and dual variables was used to prompt the global convergence of the method, where the penalty parameter is updated adaptively. Without assuming any regularity of constraints such as feasibility or constraint qualification, the method is proved to be of strong global convergence. Preliminary numerical results demonstrates that the method can not only be efficient for well-posed feasible problems, but also is applicable for some ill-posed feasible problems and some even infeasible problems.

\


\end{document}